\documentclass[11pt]{amsart}
\usepackage{latexsym,amssymb,amsmath}
\usepackage{amscd}
\usepackage[mathscr]{eucal}
\usepackage{verbatim}
\usepackage{epsfig}
\def\ctimes{\times \cdots\times}
\def\wcdots{\ww \cdots\ww}

\def\s{\sigma}

\def\ep{\epsilon}

\def\trank{\text{rank}}

\def\BC{\mathbb C}

\def\BP{\mathbb P}
\def\pp#1{\mathbb P^{#1}}

\def\fgl{\mathfrak g\mathfrak l}

\def\pp#1{{\mathbb P}^{#1}}
\def\tdim{\rm dim}
\def\hd{,...,}
\def\ww{\wedge}
\def\upperp{{}^\perp}

\def\inv{{}^{-1}}

\def\cS{{\mathfrak S}}

\def\11{\mathbf 1}
\def\PP{\mathbb P}

\def\FS{{\mathfrak S}}

\def\l{\lambda}
\def\a{\alpha}

\def\b{\beta}
\def\g{\gamma}
\def\s{\sigma}

\def\d{\delta}

\def\ot{{\mathord{\,\otimes }\,}}
\def\op{{\mathord{\,\oplus }\,}}
\def\otc{{\mathord{\otimes\cdots\otimes}\;}}
\def\cwedge{{\mathord{\wedge\cdots\wedge}\;}}
\def\ctimes{{\mathord{\times\cdots\times}\;}}

\def\lra{{\mathord{\;\longrightarrow\;}}}
\def\ra{{\mathord{\;\rightarrow\;}}}

\def\La#1{\Lambda^{#1}}

\def\tdim{\text{dim}\,}

\def\tmin{\text{ min }}

\def\thom{\text{ Hom }}
\def\trank{\text{rank}\,}

\def\tdet{\text{det}}

\newtheorem{theo}{Theorem}

\newtheorem{lemm}[theo]{Lemma}

\newtheorem{theorem}{Theorem}[section]
\newtheorem{proposition}[theorem]{Proposition}

\newtheorem{corollary}[theorem]{Corollary}

\theoremstyle{definition}
\newtheorem{definition}[theo]{Definition}

\newtheorem{example}[theo]{Example}

\theoremstyle{remark}
\newtheorem{remark}[theorem]{Remark}

\begin{document}

\title{Generalizations
of Strassen's equations for secant varieties of Segre varieties}
\author{J.M. Landsberg 
 \and L. Manivel}
%\date{December 2005}
\begin{abstract} We 
define many new examples  of modules of equations
for secant varieties of
Segre varieties that generalize Strassen's commutation equations
\cite{Strassen}.
Our modules of equations are obtained by constructing
subspaces of matrices from tensors that satisfy various
commutation properties.
\end{abstract}
\thanks{Landsberg supported by NSF grant DMS-0305829}

\maketitle
 
\section{Introduction}

Let $V,A_1\hd A_n$ be vector spaces over an algebraically
closed field $K$ of characteristic zero, and
let
$$Seg(\BP A_1\ctimes \BP A_n)\subset \BP (A_1\otc A_n)$$ denote
the {\it Segre variety} of decomposable tensors
inside $\BP (A_1\otc A_n)$.

Let $X \subset\BP V$ be a projective variety. Define $\s_r=\s_r(X)$,
 the {\it variety of
secant $\pp{r-1}$'s to $X$} by
$$
\s_r(X)=\overline{ \cup_{x_1\hd x_r\in X}\BP_{x_1\hd x_r }}
$$
where $\BP_{x_1\hd x_r }\subset\BP V$ denotes the linear space
spanned by $x_1\hd x_r$ (usually a $\BP^{r-1}$).

\smallskip 
 
For  applications
to computational complexity, algebraic statistics
and other areas, one would like to have defining
equations for secant varieties of triple Segre products,
in particular because
  the border rank $r$ of a bilinear map
$T: A^*\times B^*\ra C$ is the smallest $r$
such that $[T]\in \s_r(Seg(\BP A\times \BP B\times \BP C)$.
Here, and throughout this paper, \lq\lq defining equations\rq\rq\
refers to set theoretic
defining equations.

\smallskip

Defining equations are known only for the following cases:
all secant varieties of the two factor Segre (classical: these are
just the $(r+1)\times (r+1)$ minors of the space of $a\times b$ matrices),
 the $n$-factor Segre itself
$\BP A_1\ctimes \BP A_n\subset \BP (A_1\otc A_n)$ 
(classical), its first secant variety 
$\s_2(\BP A_1\ctimes \BP A_n)$  \cite{LMsec}, 
$\s_3(\pp {a-1}\times \pp {b-1}\times \pp {c-1})$ \cite{LWsec},
$\s_4(\pp 2\times \pp 2\times \pp 2)$ \cite{Strassen},
$\s_r(\pp 1\times \pp{b-1}\times \pp{c-1})$ \cite{LWsec}
and several cases of 
the last nontrivial secant variety
of $\pp 2\times \pp{b-1}\times\pp{b-1}$ when the last nontrivial secant
variety is a hypersurface  \cite{Strassen}.

Segre products and their secant varieties are invariant
under the action of the group $G=GL(A_1)\ctimes GL(A_n)$ and
thus their defining equations are best described as $G$-modules.
In \cite{LMsec} we explained how one can systematically find
$G$-modules in the ideal of the secant varieties using representation
theory. We also observed that   the expressions
even for highest weight vectors in the modules become too complicated
to write down explicitly very quickly, so there are severe limits
to the systematic approach.

\smallskip

The equations for $\s_2(\BP A_1\ctimes \BP A_n)$ may be thought of as those coming
from   the two factor case, that is, as minors
of ordinary matrices by considering, e.g.,
$A\ot B\ot C$ as $A\ot (B\ot C)$ and
taking the minors of the resulting $a\times bc$
matrix and permutations of such.

Strassen defined equations for $\s_3(\BP A\times \BP B\times \BP C)$
when $b=c$ and $a=3$ by choosing a basis of $A^*$ and contracting
tensors to obtain subspaces of $B\ot C$, and finding closed
conditions on such subspaces coming from tensors in $\s_3(\BP A\times \BP B\times \BP C)$.
From this perspective, one could look for other closed
conditions on such subspaces, which is one way to view our generalizations.

\smallskip

Another perspective on the equations for secant varieties
 is that
in general, if $X\subset Y$, then $\s_r(X)\subseteq \s_r(Y)$
and the equations for
$\s_2(\BP A\times \BP B\times \BP C)$ comes from  the observation that $Seg(\BP A\times \BP B\times \BP C)\subset
Seg(\BP A\times \BP (B\ot C))$.

More generally, one should look for natural
varieties, whose defining equations are easily described, that
contain $\s_r(\BP A\times \BP B\times \BP C)$.
From this perspective, our new equations are induced by
equations of  various
types of varieties of subspaces of matrices that satisfy
certain commutation properties.

For a partition $\pi$ of $d$, we let
$S_{\pi}A$ denote the corresponding irreducible
$GL(A)$ module and $\Lambda_{\pi}A=S_{\pi'}A$ where
$\pi'$ is the conjugate partition to $\pi$.
Our main result, theorem \ref{thm31}, may be phrased
as follows:

\smallskip

  {\it For each $r$, and $s$ sufficiently
small ($s\leq r/2$ if $r$ is even, $s\leq r/3$ if $r$ is odd),
we describe an explicit realization of the module
$$S_{r-s,s,s}A\ot \Lambda_{r,s}B\ot \Lambda_{r,s}C
\subset S^{r+s}(A\ot B\ot C)
$$
as a module of equations of $\s_r(\BP A^*\ot \BP B^*\ot \BP C^*)$,
and   each of these modules is independent in the ideal of $\s_r(\BP A^*\ot \BP B^*\ot \BP C^*)$.}

\smallskip

(We often reverse the roles of vector spaces with their
duals to eliminate $*$-s from the modules defining equations.)
 
The determination of the
generators of the ideal of $\s_3(\pp{a-1}\times \pp{b-1}\times \pp{c-1})$
in \cite{LWsec} relies on a computer calculation to prove the
$\s_3(\pp 2\times \pp 2\times \pp 2)$ case is generated by 
Strassen's equations, and this computer calculation was originally
announced in \cite{GSS}. In \S\ref{s3s4sect} we give a computer
free proof that the modules inherited from Strassen's equations
give set-theoretic defining equations for
$\s_3(\pp{a-1}\times \pp{b-1}\times \pp{c-1})$.
A key point in our proof is the irreducibility of the 
variety of pairs of commuting matrices. This irreducibility fails for triples
of commuting matrices. The following natural question appears to
be closely related to our problem: {\it Find equations
that characterize the irreducible component
of the variety of triples of commuting matrices containing triples of
regular semisimple matrices as an open subset.}

\medskip

One can put our investigation in the broader context
of the study of the geometry
of orbit closures: let $G$ be a complex semi-simple group, let
$V=V_{\l}$ be an irreducible $G$ module of highest weight $\l$.
Then Kostant showed that the ideal of the closed orbit
$G.[v_{\l}]=G/P\subset \BP V$ is generated in degree two
by $V_{2\l}\upperp\subset S^2V^*$. If we consider other
$G$-varieties in $\BP V$, what can we say about their
defining equations?

\subsection{Overview}
In \S\ref{inhersect} we review inheritance and remark that using
{\it subspace varieties} (defined in the section) the problem
of determining defining equations of secant varieties
of Segre varieties is reduced to the case of
$\s_r(\pp{r-1}\ctimes \pp{r-1})$.
In \S\ref{Strassensect} we review Strassen's equations for
$\s_r(\pp 2\times \pp {b-1}\times \pp{b-1})$, reformulate
them more invariantly, and describe the modules of equations
associated to his conditions. In \S\ref{rssect} we generalize
Strassen's equations and state our main result,
theorem \ref{thm31}. In \S\ref{s3s4sect} we show that
our generalizations  significantly reduce the problem of
determining defining equations in some cases, in particular
solving it when $r=3$. In \S\ref{coercesect} we generalize
our approach further and put it in a larger context, that
of a class of contractions we call {\it coercive}. Finally
in \S\ref{nontrivsect} we finish the proof of   
  theorem \ref{thm31}, showing that many of the new 
modules of equations
we defined are indeed nontrivial.

\section{Inheritance and subspace varieties}\label{inhersect}
\subsection{Inheritance}\label{inheritsubsect}
We review some facts from \cite{LMsec}.
The varieties $ \s_r(\BP A_1^*\times \cdots \times \BP A_n^*)$,
  are    invariant under the action of the group $G= 
GL(A_1)\ctimes GL(A_n)$. Thus their ideals are given by direct sums
of irreducible submodules
  $S_{\pi_1}A_1\otc S_{\pi_n}A_n\subset S^d(A_1\ot \cdots \ot A_n)$,
where each $\pi_j$ is a partition of $d$.  
 If $\tdim A_j=a_j$ then $\pi_j$ can have at most $a_j$
parts. We let $l(\pi)$ denote the number of parts
of the partition $\pi$. 
For a variety $X\subset \BP V$, we let 
$I_d(X)\subset S^dV^*$ denote the component of the  ideal
of $X$ in degree $d$.

\begin{proposition}\label{inheritprop}\cite{LMsec} If an
irreducible module $S_{\mu_1}A_1\otc S_{\mu_n}A_n\subset
I_d(\s_r(\BP A_1^*\ctimes \BP A_n^*))$, then for all vector spaces $A_j'\supseteq A_j^*$, we
have $(S_{\mu_1}A_1'\ot \cdots\ot  S_{\mu_n}A_n')^*\subset
I_d(\s_r(\BP A_1'\times \cdots \times  \BP A_n'))$.

Moreover, a module $(S_{\mu_1}A_1'\ot \cdots \ot S_{\mu_n}A_n')^*$
where the length of each $\mu_j$ is at most $a_j$ is in
$I_d(\s_r(\BP A_1'\times \cdots \times  \BP A_n'))$ iff the
corresponding module is in $I_d(\s_r(\BP A_1 \times \cdots \times
\BP A_n ))$.
\end{proposition}

 Thus a
copy of a module $S_{\mu_1}A_1\ot \cdots \ot S_{\mu_n}A_n$ will be
in $I(\s_r(\pp {r-1}\times \cdots \times \pp{r-1}))$
 iff the corresponding copy of the module
$S_{\mu_1}\BC^{l(\mu_1)}\ot \cdots \ot S_{\mu_n}\BC^{l(\mu_n)}$ is
in the ideal of $\s_r(\pp {l(\mu_1)-1}\times \cdots \times
\pp{l(\mu_n)-1})$.

\subsection{Subspace varieties}
Let $Sub_{b_1\hd b_n}\subset \BP (A_1^*\ot\cdots \ot A_n^*)$
denote
  the set of tensors $T$ such that there
exists subspaces $B_j \subseteq A_j^*$ with $\tdim B_j=b_j$ and
$T\in B_1\ot \cdots \ot B_n$. $Sub_{b_1\hd b_n}$ is Zariski closed
  and     its ideal is   easy to describe.
$I_d(Sub_{b_1\hd b_n})$ is the direct sum of the modules
$S_{\mu_1}A_1\otc S_{\mu_n}A_n$ such that $S_{\mu_1}A_1\otc
  S_{\mu_n}A_n\subset S^d(A_1\ot\cdots\ot A_n)$ and the
length of some $\mu_j$ is greater than $b_j$.
(In \cite{LWsec} we prove the generators of the ideal
are indeed the expected ones.)

Assuming all the $b_j$ are equal to say $b_0$, then
$Sub_{b_0\hd b_0}$ is defined by equations of degree
$b_0+1$, namely all the modules in $S^{b_0+1}(A_1\otc A_n)$
containing an exterior power of some $A_j$.
In other words, as a set, $Sub_{r\hd r}$ is the intersection
of all the $r$-th secant varieties of  flattenings of the form
$A_i\ot (A_1\otc \hat A_i\otc A_n)$.

In particular,  $\s_r(\BP A_1^*\times\cdots \times \BP A_n^*)\subset
Sub_{b_1\hd b_n}$ for all $b_1\hd b_n$ with $b_i\geq r$.  
We summarize the above discussion:

\begin{proposition}\label{subspaceprop} Defining equations for $\s_r(\BP
A_1^*\times\cdots\times \BP A_n^*)$, when $\tdim A_j^*\geq r$  
may be obtained from the union of the the modules inherited from 
defining equations for  $\s_r(\pp {r-1}\times \cdots \times \pp{r-1})$
and defining equations for  $Sub_{r\hd r}$. 
\end{proposition}

\begin{remark} For ordinary matrices,
i.e., points in the tensor product
of two vector spaces,  there is just one notion
of rank, but it has several generalizations to
tensor products of several vector spaces.
The first is the minimum number of monomials required to
express a given tensor
as a sum of monomials, which is now commonly
called the {\it rank} of the tensor. The second
is the smallest secant variety of the Segre variety
in which the tensor lies, which is called the
{\it border rank} of the tensor. A third notion
comes from {\it Cayley's hyperdeterminant}, a higher dimensional
generalization of the determinant. Already for $\pp 1\times\pp 1\times\pp 1$
this notion diverges from the previous two, in the sense that
every tensor in $\pp 1\times\pp 1\times\pp 1$ has border rank at
most two, but the zero set of the hyperdeterminant is a quartic hypersurface.
For $\pp 2\times\pp 2\times \pp 2$ the hyperdeterminant describes
an irreducible hypersurface of degree $36$ whereas $\s_4(\pp 2\times \pp 2\times \pp 2)$ is a hypersurface of degee 9. The hyperdeterminants induce
\lq\lq hyper-minors\rq\rq\  
by inheriting the corresponding modules, but the zero sets of these appear to have
little relation with secant varieties.
A fourth notion generalizes to the {\it subspace varieties}, because
  $T\in A\ot B$ has rank $r$ iff there exist
$A'\subset A$, $B'\subset B$, both of dimension 
$r$, with $T\in A'\ot B'$. Tensors of border rank $r$ are in general only contained
in $Sub_{r\hd r}$.

\end{remark}

\section{Strassen's equations}\label{Strassensect}

\subsection{Strassen's theorem}
For a tensor $T\in A\ot B\ot C$ and $\a\in A^*$, let 
$T_{\a} \in B\ot C$ denote the contraction of $T$ with $\a$.

\begin{theorem}[Strassen]\cite{Strassen} Let $3\leq a\leq b=c\leq r$. Let
$T\in A\ot B\ot C$ and $\a\in A^*$ be such that $\trank T_{\a}=b$. For all $\a^1,\a^2\in A^*$,
consider the linear maps $T_{\a,\a^j}: B\ra B$ by considering
$T_{\a}: C^*\ra B$ and $T_{\a,\a^j}= T_{\a^j}T_{\a}\inv$. 
If $[T]\in \s_r(\BP A\times \BP B\times \BP C)$, then
$$
\trank [T_{\a,\a^1},T_{\a,\a^2}]\leq 2(r-b).
$$
Moreover 
for a generic tensor $T\in A\ot B\ot C$, $[T_{\a,\a^1},T_{\a,\a^2}]$
is of maximal rank.
\end{theorem}

This theorem (together with an easy application of
Terracini's lemma) implies $\s_4(\pp 2\times \pp 2\times \pp 2)$ is  
a hypersurface.
It also implies that the
border rank of the multiplication of $m\times m$
matrices is at least $\frac {3m^2}2$.
Here is a proof that is essentially Strassen's, rephrased
more invariantly to enable generalizations.

\begin{proof}
 First note that it is sufficient to prove the result for
$T$
of the form
$T=a_1b_1c_1+\cdots +a_rb_rc_r$
 as these form a 
Zariski open subset of the irreducible variety
$\s_r$. Here $a_j\in A$ etc... and $a_jb_jc_j=a_j\ot b_j\ot c_j$.
Fix an auxiliary vector space $D\simeq \BC^r$ and write
$T_{\a}: C^*\ra B$ as a composition of maps
$$
\begin{CD}
C^* &@>i>>& D &@>\d_{\a}>>& D &@>p>>& B.
\end{CD}
$$

 To see this explicitly, if $T=a_1b_1c_1+\cdots +a_rb_rc_r$ and
we assume $b_1\hd b_b$, $c_1\hd c_b$ are bases of $B,C$, then
letting $d_1\hd d_r$ be a basis of $D$, we have
$i(\eta)=\sum_{j=1}^r\eta(c_j)d_j$, $\d_{\a}(d_j)=\a(a_j)d_j$,   for $1\leq s\leq b$ we have
$p(d_s)=b_s$, and for $b+1\leq x\leq r$, writing $b_x=\xi^s_xb_s$,
then  we have $p(d_x)=\xi^s_x b_s$.

Let $D'=i(C^*)$,
write $i': C^*\ra D'$ and set $p_{\a}:=p\mid_{\d_{\a}(D')}$, so
$p_{\a}:\d_{\a}(D')\ra B$ is a linear isomorphism.
Then we may write $T_{\a}\inv = (i')\inv \d_{\a}\inv p_{\a}\inv$.

Note that $\trank [T_{\a,\a^1},T_{\a,\a^2}]=\trank(
T_{\a^1} T_{\a}\inv T_{\a^2}-T_{\a^2}T_{\a}\inv T_{\a^1})$
because $T_{\a}$ is invertible.
We have
\begin{align*}
T_{\a^1}& T_{\a}\inv  T_{\a^2}-T_{\a^2}T_{\a}\inv T_{\a^1}\\
&=
(p\d_{\a^1}i')((i')\inv \d_{\a}\inv p_{\a}\inv)(p\d_{\a^2}i')
-
(p\d_{\a^2}i')((i')\inv \d_{\a}\inv p_{\a}\inv)(p\d_{\a^1}i')\\
&=
p[
 \d_{\a^1} \d_{\a}\inv p_{\a}\inv p\d_{\a^2} 
-
 \d_{\a^2}  \d_{\a}\inv p_{\a}\inv p\d_{\a^1}]i'\\
&=
p\d_{\a}\inv[
 \d_{\a^1}  p_{\a}\inv p\d_{\a^2} 
-
 \d_{\a^2}   p_{\a}\inv p\d_{\a^1}]i'
\end{align*}
where the last equality holds because the $\d_{\a}$'s commute.

Now $p_{\a}\inv p\mid_{\d_{\a}(D')}=Id$, so write
$D= \d_{\a}(D')\op D''$, where we choose any complement
to $\d_{\a}(D')$ in $D$. We have $\tdim D''= r-b$ and 
we may write $p_{\a}\inv p=Id_{\d_{\a}(D')} +f$
for some map $f: D''\ra D$. Thus
$$T_{\a^1} T_{\a}\inv T_{\a^2}-T_{\a^2}T_{\a}\inv T_{\a^1}
=
p\d_{\a}\inv[
 \d_{\a^1}  f\d_{\a^2} 
-
 \d_{\a^2}   f\d_{\a^1}]i'
$$
and is therefore of rank at most $2(r-b)$.\end{proof}

\subsection{Towards a more invariant formulation of Strassen's theorem}
As stated, there are several undesirable aspects to
Strassen's equations:
 the choices of $\a,\a^1,\a^2$, 
the requirement that $\a$ is such that $T_{\a}$ invertible,
and the way the equations are written makes it difficult to
see what equations will be inherited from them when
we increase the dimensions of the spaces.
Moreover, say $a=b=c$, then we can clearly change
the roles of the spaces - are the new equations so
obtained redundant or not?

A first step towards resolving these issues
is to reconsider matrix multiplication and
inverses more invariantly.

\smallskip

For a linear map $f: V\ra W$, let $f^{\ww k}: \La k V\ra \La k W$ denote
the induced linear map.
If $\tdim V=\tdim W=n$ and $f$ is invertible, then
as a tensor $f^{\ww n-1}=(f\inv)^t \ot \tdet (f)$.  Recall   that for a vector space of dimension
$n$, that $\La{n-1}V=V^*\ot \La n V$, so
if $V,W$ have dimension $n$, then $\La{n-1}V^*\ot \La{n-1}W=
\thom (W,V)\ot \La b V\ot \La b W$. 
$f^{\ww n-1}$ has the advantage over $f\inv$ of being defined
even if $f$ is not invertible.

For $T\in A\ot B\ot C$, let
 $T^{\a}:=T_{\a}^{\ww b-1}\in \La{b-1}B\ot \La{b-1}C=
\La{b-1}B\ot C^*\ot \La b C$.  We may contract $T^{\a} \ot T_{\a^j}
\in\La{b-1}B\ot C^*\ot \La b C\ot   B\ot C$ to
an element 
$$T^{\a}_{ \a^j}\in \La b B\ot C^*\ot \La b C\ot C
=C^*\ot C\ot \La b B\ot \La b C.
$$
Now consider 
$$T^{\a}_{ \a^1}\ot T^{\a}_{ \a^2}
\in  C^*\ot C\ot C^*\ot C\ot (\La b B)^{\ot 2}\ot (\La b C)^{\ot 2}
$$
and contract  on the second and third factors
to obtain an element of 
$C^*\ot C\ot (\La b B)^{\ot 2}\ot (\La b C)^{\ot 2}$. This contraction
of course corresponds to matrix multiplication, as does
contraction in the first and fourth factor, which corresponds
to multiplying the matrices in the opposite order.  
We do both contractions and take their difference 
and call the result
$$
[T^{\a}_{ \a^1},T^{\a}_{ \a^2}]\in C^*\ot C\ot  (\La b B)^{\ot 2}\ot (\La b C)^{\ot 2}.
$$
  Strassen's theorem states that the rank of $[T^{\a}_{ \a^1},T^{\a}_{ \a^2}]$
is at most $2(r-b)$.

Equivalent to Strassen's observation
that $\trank [T_{\a, \a^1},T_{\a, \a^2}]=
\trank(T_{\a^1}T_{\a}\inv T_{\a^2}-T_{\a^2}T_{\a}\inv T_{\a^1})$, 
we can get away with
a lower degree tensor by just contracting 
once with $T^{\a}$ to get elements of $B\ot C\ot \La b B\ot \La b C$.

To eliminate the choices of $\a,\a_1,\a_2$,
we may consider the tensor $T^{\a} $
without having chosen $\a$ as $T^{(\cdot)} \in S^{b-1}A\ot \La{b-1}B\ot \La{b-1}C$,
which is obtained as the projection of $(A\ot B\ot C)^{\ot b-1}$ 
to the subspace $S^{b-1}A\ot \La{b-1}B\ot \La{b-1}C$. Similarly 
$T_{\a^j}\in B\ot C$, may be thought of as
$T_{(\cdot)}\in A\ot B\ot C$.
We then contract 
$$T^{(\cdot)} \ot T_{(\cdot)}\ot T_{(\cdot)}
\in \La{b-1}B\ot \La{b-1}C\ot B\ot C\ot B\ot C \ot (S^{b-1}A\ot A\ot A)
$$
in two different ways, first contracting the first factor
with the third and the second with the sixth to obtain
an element of $B\ot C\ot (S^{b-1}A\ot A\ot A)\ot \La b B\ot \La b C$, then contracting the first with
the fifth and the second with the fourth. We then
take the difference of the two to obtain
an element of $B\ot C\ot \La b B\ot \La bC\ot (S^{b-1}A\ot A\ot A)$. Call the resulting tensor
$\phi(T)$, i.e.,
$$
\phi\in  (A^*\ot B^*\ot C^*)^{\ot b+1}\ot  ( S^{b-1}A\ot \La 2 A) \ot   \La b B\ot B\ot 
 \La b C\ot C
$$
and in fact descends to be an element of
$S^{b+1}(A^*\ot B^*\ot C^*)\ot  ( S^{b-1}A\ot \La 2 A) \ot   \La b B\ot B\ot 
 \La b C\ot C$.

The proof of Strassen's theorem may be rephrased in this
language. We leave this as an entertaining exercise for the
reader. (Hint: the Pl\"ucker relations for the Grassmannian
$G(2,r)$ furnish the key to showing the bound on
the rank of the commutator.)  

\subsection{Strassen's equations as modules}
We  first determine
which modules in
$$
\La 2 A\ot S^{b-1}A\ot \La bB\ot B\ot C\ot \La bC
$$ 
map nontrivially
into $S^{b+1}(A\ot B\ot C)$, when we
use $\phi$ to compose the
inclusion 
$$\La 2 A\ot S^{b-1}A\ot \La bB\ot B\ot C\ot \La bC
\subset (A\ot B\ot C)^{\ot b+1}$$
 with the projection
$(A\ot B\ot C)^{\ot b+1}\ra S^{b+1}(A\ot B\ot C)$.

Since here $b=\tdim B=\tdim C$, we have
$$\La 2 A\ot S^{b-1}A\ot \La bB\ot B\ot C\ot \La bC
=(S_{b,1}A\op S_{b-1,1,1}A)\ot  \Lambda_{b,1}B \ot \Lambda_{b,1}C
$$
so there are two possible modules. By \cite{LWsec},
$S_{b,1}A\ot \Lambda_{b,1}B\ot \Lambda_{b,1}C$ does
not occur in the ideal of
$\s_r(\pp 1\times\pp{b-1}\times\pp{c-1})$ so by inheritance
it cannot occur when $\tdim A>2$ either, so we are reduced
to a unique module.

Taking minors corresponds to taking exterior powers in the $B,C$
factors and we conclude:

\begin{proposition}
As modules, the equations that imply
$$
\trank [T^{\a}_{\a^1},T^{\a}_{\a^2}]\leq 2(r-b).
$$
for all choices of $\a,\a_1,\a_2\in A^*$
correspond to the image of the inclusion 
via $\phi$ of 
$$S^{2(r-b)+1}(S_{b-1,1,1}A)\ot \La {2(r-b)+1}( \Lambda_{b,1}B)\ot 
\La {2(r-b)+1}(\Lambda_{b,1}C)
$$
into $S^{(2(r-b)+1)(b+1)}(A\ot B\ot C)$.
When $r=b$ we obtain the single module
$$S_{b-1,1,1}A\ot \Lambda_{b,1}B\ot \Lambda_{b,1}C.$$
\end{proposition}

For $\pp 2\times \pp 2\times \pp 2$ we obtain
the same modules regardless of which factor we use to make
the projections - all are 
the same copy of $S_{211}A\ot S_{211}B\ot S_{211}C$ for
the case of $\s_3$ and of $S_{333}A\ot S_{333}B\ot S_{333}C$ for
the case of $\s_4$. This
redundancy fails for larger dimensional projective
spaces.

\subsection{Example of Strassen's equations}
We write down a basis of the modules of polynomials corresponding to    $S_{b-1,1,1}A\ot \Lambda_{b,1}B\ot \Lambda_{b,1}C$. Let $\a_i,\a_j,\a_k\in A^*$, let  $\b_1\hd \b_b$, $\xi_1\hd \xi_b$ be bases of
$B^*,C^*$. Consider the tensor
$$
P^{i,j|k}_{s,t}=\a_i\ww\a_j\ot (\a_k)^{b-1}\ot \b_1\wcdots \b_b
\ot \b_s\ot \xi_1\wcdots \xi_b\ot \xi_t
$$
Applying $\phi$ we obtain (ignoring scalars)
\begin{align*}
&(\a_2\ot \a_3-\a_3\ot \a_2)\ot (\a_1)^{b-1}\ot
(\sum_j (-1)^{j+1}\b_{\hat \j} \ot \b_j\ot \b_s)\ot
(\sum_k (-1)^{k+1}\xi_{\hat k} \ot \xi_k\ot \b_t)\\
&=
(-1)^{j+k}[
((\a_1)^{b-1}\ot \b_{\hat \j}\ot \xi_{\hat k} )
\ot (\a_2\ot \b_j\ot \xi_t)\ot (\a_3\ot \b_s\ot \xi_k)\\
&\ \ \ -
((\a_1)^{b-1}\ot \b_{\hat \j}\ot \xi_{\hat k} )
\ot (\a_3\ot \b_j\ot \xi_t)\ot (\a_2\ot \b_s\ot \xi_k)].
\end{align*}
If we choose dual bases for $A,B,C$ and write
$$T=\sum_l a_l\ot X_l $$
 where
the $a_l$ are dual to the $\a_l$ and $X_l$
are represented as
$b\times b$ matrices with respect to the dual bases of $B,C$,
then
$$
P^{i,j|k}_{s,t}(T)=
\sum_{u,v}(-1)^{u+v}(\tdet X^{\hat u}_{k,\hat v})
(X^u_{i,t}X^s_{j,v}-X^s_{i,v}X^u_{j,t})
$$
where $X^{ \hat u}_{j,\hat v}$ is $X_j$ with its $u$-th row and
$v$-th column removed.

\section{Generalizations of Strassen's conditions}\label{rssect}
We now generalize Strassen's equations using our new perspective.
Recall that the key point
for Strassen's equations
was that contracting 
a tensor $T\in A\ot B\ot C$ in two different ways yielded tensors that  
almost
commute  when $T\in \s_r$.

Consider, for $s,t$ such that $s+t \leq b$ and $\a,\a_j\in A^*$,
the tensors
$$T_{\a_j}^{\ww s}\in \La sB\ot \La s C,\ 
T_{\a}^{\ww t}\in \La tB\ot \La t C
$$
(our old case was $s=1, t=b-1$).
We may contract $T_{\a}^{\ww t}\ot T_{\a_1}^{\ww s}\ot T_{\a_2}^{\ww s}$ to
obtain   elements of $\La {s+t}B\ot \La {s+t}C\ot \La s  B\ot \La s C$ in two
different ways, call these contractions $\psi^{s,t}_{\a,\a_1,\a_2}(T)$ and
$\psi^{s,t}_{\a,\a_2,\a_1}(T)$. 

Now say   we may write $T=a_1\ot b_1\ot c_1+\cdots + a_r\ot b_r\ot c_r$
for elements $a_i\in A$, $b_i\in B$, $c_i\in C$.
We have
$$\psi^{s,t}_{\a,\a_1,\a_2}(T)=\sum_{|I|=s,|J|=t,|K|=s}
\langle a_{I},\a_1\rangle\langle a_{J},\a\rangle \langle
a_{K},\a_2\rangle (b_{I+J}\ot b_{K})\ot (c_{I}\ot
c_{J+K}),$$ 
where we used the notation $a_{I+J}=a_I\wedge a_J$ etc.. For
this to be nonzero, we need $I$ and $J$ to be disjoint subsets of
$\{1,\ldots ,r\}$. Similarly, $J$ and $K$ must be disjoint.
If $s+t=r$ this implies $I=K$.
We conclude:

\begin{proposition} For $T\in \s_{s+t}(\BP A\times \BP B\times\BP C)$, 
for all $\a,\a_1,\a_2\in A^*$
$$
\psi^{s,t}_{\a,\a^1,\a^2}(T)-\psi^{s,t}_{\a,\a^2,\a^1}(T)=0.
$$ 
\end{proposition}

We  have the bilinear map
$$
(\La 2(S^sA)\ot S^tA)^*\times (A\ot B\ot C)^{\ot 2s+t}\ra \La{s+t}B\ot \La{s+t}C
\ot \La sB\ot \La sC.
$$
whose image is $\psi^{s,t}_{\a,\a^1,\a^2}(T)-\psi^{s,t}_{\a,\a^2,\a^1}(T)$.
We rewrite it as a polynomial map
$$
\Psi^{s,t} :  A\ot B\ot C \ra 
(\La 2(S^sA)\ot S^tA)\ot 
\La{s+t}B\ot \La{s+t}C
\ot \La sB\ot \La sC.
$$
If we want to consider polynomial equations on $A\ot B\ot C$,
they are the image of the transpose of $\Psi^{s,t}$.
So just as with Strassen's equations, we no longer need to make
choices of elements of $A^*$.
The only catch is we don't yet know whether or not  
$\Psi^{s,t}(T)$ is identically
zero for all tensors $T$. This is addressed in \S\ref{nontrivsect}.

We write $r=s+t$ and call the image of $\Psi^{s,r-s}$ the
{\it $(r,s)$-coercive equations}.

\smallskip

The modules for
the $(r,s)$-coercive equations
are the irreducible submodules of
$$\Lambda^2(S^sA)\ot S^{r-s}A \ot\La{r}B\ot\La{s}B
\ot\La{s}C\ot\La{r}C$$
 that map isomorphically into
$S^{r+s}(A\ot B\ot C)$ under the transpose of $\Psi^{s,t}$.
There are many such submodules and we can describe
them explicitly (see the formulas \eqref{LRrules} below), but
  there is no easy to implement formula for the decomposition
of $S^{r+s}(A\ot B\ot C)$. There are certain
modules that  are easily seen to occur in both $S^{r+s}(A\ot B\ot C)$ and
$\Lambda^2(S^sA)\ot S^{r-s}A \ot\La{r}B\ot\La{s}B
\ot\La{s}C\ot\La{r}C$, and we will show
that at least most of the time, these modules map isomorphically,
so that most of the $(r,s)$-coercive conditions lead to  
  nontrivial equations.

\begin{theorem}\label{thm31}
For $s$ odd, $r$ even,  and $2s\leq r$,
or $r,s$ odd and $3s\leq r$,  the multiplicity one component of
$S_{r+s}(A\ot B\ot C) $ of type $S_{r-s,s,s}A 
\ot\Lambda_{r,s}B  \ot\Lambda_{r,s}C $
induced from the $(r,s)$-coercive equations is a 
nontrivial set of equations of
$\s_r(\PP A^*\times\PP B^*\times\PP C^*)$.
All these modules of equations for $\s_r$ are independent  
elements of the ideal of $\s_r$ as $s$ varies.
\end{theorem}
 
The proof of nontriviality is given in \S\ref{nontrivsect}.
To see the independence, consider
equations of degrees $r+s_1$ and $r+s_2$ with $s_1<s_2$.
Were the second set induced by the first,
the corresponding tableau  for
$S_{r-s_1,s_1,s_1}A$ would have to fit inside
the tableau for $S_{r-s_2,s_2,s_2}A$.
But since $r-s_1>r-s_2$ the first
tableau has a longer
first row than the second.

\section{$Comm^r_A$}\label{s3s4sect}

We   study the special case $s=1$. Then
we have the decompositions
\begin{align*}
 \La 2 A\ot S^{r-1}A&= S_{r,1}A\op S_{r-1,1,1}A,\\
B\ot \La{r}B&=\La{r+1}B\op \Lambda_{r,1}B,\\
C\ot \La{r}C&=\La{r+1}C\op \Lambda_{r,1}C.
\end{align*}
We may ignore the modules containing an $(r+1)$-st exterior
power as we already know all those are contained in
the ideal, and  we may eliminate
 $S_{r,1}A\ot \Lambda_{r,1}B\ot \Lambda_{r,1}C$ as above.  
Thus we are reduced to studying
the modules inherited from Strassen's equations. 

\begin{definition}We let $Comm^r_A\subset \BP (A\ot B\ot C)$ be
the set of tensors $T$ such that $\Psi^{1,r-1}(T)=0$. This set
is Zariski closed  whose ideal
is generated by the image of the transpose of $\Psi^{1,r-1}$. The case
$a=3$,
$r=b=c$ corresponds to the tensors obeying Strassen's
commutation condition $[T_{\a,\a_1},T_{\a,\a_2}]=0$ for
all $\a,\a_1,\a_2\in A^*$ such that $T_{\a}$ is invertible.
\end{definition}

Note   that these equations are of the minimal
degree $r+1$ (see \cite{LMsec}).

\begin{proposition} \label{a=3}
For $a=3\leq r\leq b,c$, 
$$Comm^r_{A}=\s_r(\BP A^*\times \BP B^*\times
\BP C^*)\cup Sub_{3,r-1,r}\cup Sub_{3,r ,r-1}.
$$
\end{proposition}

\begin{proof}
The set of defining equations for $Comm^r_{A}$ is  
$S_{211}A\ot \Lambda_{r,1}B\ot \Lambda_{r,1}C$. In particular they all
involve terms containing partitions of length $r$ in $B$ and $C$,
thus they vanish on $Sub_{3,r,r-1}\cup Sub_{3,r-1,r}$.
We also already saw that $Comm^r_{A}\supseteq \s_r$, so we
have
$$Comm^r_{A}\supseteq \s_r(\BP A^*\times \BP B^*\times \BP C^*)\cup Sub_{3,r,r-1}\cup Sub_{3,r-1,r}.
$$

Let $T\in Comm^r_{A}$
 be such that $T\notin  Sub_{3,r,r-1}\cup Sub_{3,r-1,r}$.
Let $B'\subset B$, $C'\subset C$ be the smallest subspaces
such that $T\in A\ot B'\ot C'$.
$T\in Comm^r_{A}$
implies that $B',C'$ both have dimension at most $r$
and $T\notin  Sub_{3,r,r-1}\cup Sub_{3,r-1,r}$,
implies further that both have dimension exactly $r$.

 Fix $\a_0\in A$
and consider the abelian subalgebra $\{
T_{\a_0, \a_1},T_{\a_0,\a_2} \}\subset End (C')$. Now the crucial point  is that
any pair of commuting matrices can be approximated by
simultaneously diagonalizable matrices. 
(This statement is the only place where we use the hypothesis
that $a=3$. It is a slightly more
precise statement than the well-known irreducibility of the
commuting variety \cite{mz}.  Note that the corresponding
statement is not true for three or more commuting matrices.) That
is, our tensor $T$ is in the closure of the set of those $T'$'s
for which we can find a basis $b_1,\ldots ,b_r$ of $B'$, and a
basis $c_1,\ldots ,c_r$ of $C'$, such that any $T'(\a)$ is a
linear combination of $b_1\ot c_1,\ldots ,b_r\ot c_r$. 
 But then we
can find $a_1,\ldots ,a_r$ in $A$, such that $T'=a_1\ot b_1\ot
c_1+\cdots +a_r\ot b_r\ot c_r$. In particular such a $T'$ belongs
to $\s_r$, hence so does $T$.
\end{proof}

Now, since $\s_3(\pp 2\times \pp 2\times \pp 1)$ is the entire
ambient space, $Sub_{3,3,2}\cup Sub_{3,2,3}\subset \s_3(\pp{a-1}\times
\pp{b-1}\times\pp{c-1})$ and we conclude:

\begin{corollary}\label{s3thm} As sets, for $a,b,c\geq 3$,
$$\s_3(\pp {a-1}\times \pp {b-1}\times \pp {c-1})=Comm^3_{A}\cap Sub_{333}
=Comm^3_{B}\cap Sub_{333}=Comm^3_{C}\cap Sub_{333}.
$$
 That is $\s_3(\pp {a-1}\times \pp {b-1}\times \pp {c-1})$
is the zero set of 
$S_{211}A\ot S_{211}B\ot S_{211}C\subset S^4(A\ot B\ot C)$
and modules in degree four containing a fourth exterior power (i.e.,
$\La 4A\ot \La 4(B\ot C)$ plus permuatations).
In particular, $\s_3$ is cut out set-theoretically by equations
of degree four.
\end{corollary}

\begin{remark} In fact, the stronger statement that
the ideal of $\s_3$ is generated by the above modules
holds, see \cite{LWsec}, but the proof relies
on a computer calculation.
\end{remark}

\begin{proposition}\label{rleq4}
For $a\leq b,c$ and $r\leq 4$, 
$$Comm^r_{A}=\s_r(\BP A^*\times
\BP B^*\times \BP C^*) \cup Sub_{a,r-1,r }\cup Sub_{a,r ,r-1}.
$$
\end{proposition}

\begin{proof}
The proof is the same as above except that at the point where
we used $a=3$ we use instead that 
  for $r\leq 4$, an $r$-dimensional abelian subalgebra
of $\fgl_r$ can be approximated by Cartan subalgebras (subalgebras
of matrices that are diagonal in some fixed basis)
\cite{IM} and we conclude
as above.
\end{proof}

\begin{remark} It is likely that $5$-dimensional
abelian subalgebra of $\fgl_5$ can be approximated by Cartan
subalgebras so proposition \ref{rleq4} should still hold for
$r=5$, \cite{IM}. On the other hand, it is not possible to approximate
$r$-dimensional abelian subalgebras of $\fgl_r$ by Cartan algebras
for $r>5$.
\end{remark}

\begin{corollary}\label{s4prop}
As sets, for $a,b,c\geq 3$,
  $\s_4(\pp {a-1}\times \pp {b-1}\times \pp {c-1})$
is the zero set of 
\begin{enumerate}
\item $(S_{311}A \ot S_{2111}B \ot S_{2111}C )
\op (S_{2111}A \ot S_{311}B \ot S_{2111}C ) 
\op (S_{2111}A \ot S_{2111}B \ot S_{311}C )
\subset S^5(A\ot B\ot C) $,
i.e., the equations of $Comm^4$.
 
\item
  equations  inherited  from $\s_4(\pp 2\times \pp 2\times \pp 3)$
 
\item
   modules  in  $S^5(A\ot B\ot C)$ containing 
a  fifth  exterior  power, i.e., the equations for $Sub_{4,4,4}$. 
\end{enumerate}
\end{corollary}

\begin{remark}  The known defining modules for
$\s_4(\pp 2\times \pp 2\times \pp 3)$ are 
$S_{321}A\ot S_{321}B\ot S_{3111}C$ in degree $6$ and
$S_{333}A\ot S_{333}B\ot S_{333}C$ in degree $9$, \cite{LMsec}.
We do not have an interpretation for $S_{321}A\ot S_{321}B\ot S_{3111}C$,
and it would be useful to have one in order to determine if
the known modules for $\s_4(\pp 2\times \pp 2\times \pp 3)$ are sufficient
to define it. 
In \cite{LMsec} there is a typographical error in the 
statement of proposition 6.3, incorrectly giving the modules
in degree six, although they are written correctly in the proof.
\end{remark}

\section{Coercive contractions}\label{coercesect}
We now place the   discussion of \S\ref{rssect} in a more general context.
Let $m$ and $k$ be integers, with $m$ even. Consider the projection
\begin{eqnarray}\nonumber
S^k(A_1\otc A_m) & \lra & \Lambda^{k}A_1\otc\Lambda^{k}A_m,
%\Lambda^{2k+1}(A_1\otc A_m) & \lra & \Lambda^{2k+1}A_1\otc
%\Lambda^{2k+1}A_m.
\end{eqnarray}
sending $T=\sum_{i}a_1^i\otc a_m^i$ to
$$\wedge^{k}T=\sum_{|I|=k}a_1^I\otc a_m^I,$$
where if $I=(i_1<\cdots <i_{k})$, then $a^I=a^{i_1}\cwedge a^{i_{k}}$.

Let  
$$T=\sum_{1\le i\le r}a_0^i\otc a_m^i\in A_0^*\otc A_m^*$$
for some vectors $a^i_j\in A_j^*$. For any $\a\in A_0$, let
$T(\a)\in A_1^*\otc A_m^*$ denote the contraction of $T$ by $\a$.
Then
$$\wedge^{k}T(\a)=\sum_{|I|=k}\langle a_0^I,\a\rangle
a_1^I\otc a_m^I.
$$  
Now consider the product
of $p$ such tensors,
\begin{align}\label{coercemap}
&\wedge^{k_1}T(\a_1)\otc \wedge^{k_p}T(\a_p) \hspace*{8cm}\hfill \\
&\nonumber = \sum_{|I_1|=k_1,\ldots , |I_p|=k_p} \langle
a_0^{I_1},\a_1\rangle\cdots \langle a_0^{I_p},\a_p\rangle
(a_1^{I_1}\otc a_1^{I_p})\otc (a_m^{I_1}\otc a_m^{I_p}).
\end{align}
Note that we put together the different terms involving wedge
powers of each $A_j^*$. 
This is because we want to take
more skew-symmetrizations, that is, we want to apply natural maps
of type 
\begin{equation}\label{contrw}\Lambda^{m_1}A_1^*\otc \Lambda^{m_t}A_1^* \ra
\Lambda^{m_1+\cdots +m_t}A_1^*
\end{equation} to our tensor.

\begin{definition}
A contraction
\begin{equation}\label{contr}\Gamma :A_0^{\times p}\times (A_0^*\otc A_m^*)
\ra (\La{k_1}A_1^*\ot \cdots \ot \La{k_1}A_m^*)\ot \cdots \ot
(\La{k_p}A_1^*\ot \cdots \ot \La{k_p}A_m^*)
\end{equation}
 given by \eqref{coercemap} followed by maps of the form \eqref{contrw}
is {\it  r-coercive}  if when restricted to 
tensors of the form $T=a^1_0\otc a^1_m+\cdots + a^r_0\otc a^r_m$ the
only nonzero terms in the right hand side of \eqref{coercemap} are
terms with $I_1=I_2$ (or more generally the only nonzero terms are
those where two of the multi-indices $I_j$  coincide).  
Generalizing our previous discussion, $r$-coercive contractions furnish equations
for $\s_r(\BP A_0^*\ctimes \BP A_m^*)$, by taking
$(\Gamma -\Gamma')(T)$ where $\Gamma'$ is the same
as $\Gamma$ only switching the roles of the coinciding multi-indices.

Such a contraction is called {\it  partially r-coercive} if 
when restricted to tensors of the form $T=a^1_0\otc a^1_m+\cdots + a^r_0\otc a^r_m$ the contracted tensor is nongeneric
 among tensors in the image of $\Gamma-\Gamma'$.
\end{definition}

Strassen's tensors $\phi$   are  partially $r$-coercive because
the contracted tensor can have rank (as a matrix) at most
$2(r-b)$ whereas a generic such matrix has rank $b$.
The tensors $\Psi^{s,t}$ are $(s+t)$-coercive and partially
$(s+t+x)$-coercive for small $x$. 

Partially coercive tensors $\Gamma$ applied to
$T\in \s_r$ for sufficiently small $r$ give rise to
tensors $\Gamma (T)$ that  
belong to some type of secant variety. Since our understanding
of higher secant varieties is quite limited in general, it
is not always clear how to use them. However, there is one case
we understand well, namely the secant varieties of two-factor
Segre varieties, which is what is used for the Strassen equations.  
 
\smallskip

Here is a more complicated example of a coercive contration:

\begin{example} Let $m=6$ and $p=7$, $k_1=k_2=k_3=r-4s$ and
$k_4=k_5=k_6=k_7=s$. Then the contraction
$\psi_{145,167,246,257,347,356}$ is $r$-coercive.
(Here the grouped indices indicate which
are to be contracted together.)  Indeed, the
contraction $145$ implies for the surviving terms that $I_1\cup
I_4\cup I_5$ is a disjoint union, in other words $I_4$ and $I_5$
are disjoint and $I_1$ is contained in the complement of their
union. Taking the other contractions into account, we see that
$I_4,I_5,I_6,I_7$ are pairwise disjoint, and that $I_1,I_2,I_3$
are contained in, hence equal to because of the cardinalities, the
complement of their union. In particular they must be equal.
\end{example}

\begin{example} Here are some further examples of partially
coercive equations, and further variants on
these should be clear to the reader.
In the propositions below we assume $b=c$. Of course
the corresponding modules induce equations when this
is not the case, but moreover when $b\leq r\leq c$
there are further modules of equations that
are induced that are not inherited.

\begin{proposition}\label{newprop}  
Let $T\in A^*\ot B^*\ot C^*$ and $\a_0,\a_1,\a,\a'\in A$. If 
$T\in \s_r(\BP A^*\times \BP B^*\times \BP C^*)$, then
$$rank\; [T^{\a_0}_{\a},T^{\a_1}_{\a'}]\le 3(r-b).$$
The relevant modules are the corresponding image of
$\La 2(S^{b-1}A)\ot \La 2 A\ot \La b B\ot B\ot \La b C\ot C$ in
$S^{2b}(A\ot B\ot C)$.
\end{proposition}

\begin{proposition}
Let $T\in A^*\ot B^*\ot C^*$ and $\a_0,\a_1\hd \a_k\in A$. If $T\in \s_r(\BP A^*\times \BP B^*\times \BP C^*)$, then for any permutation $\s\in\cS_k$,
$$rank\; \big(T^{\a_0}_{\a_1}\cdots T^{\a_0}_{\a_k}
-T^{\a_0}_{\a_{\s(1)}}\cdots T^{\a_0}_{\a_{\s(k)}}\big) \le
2(k-1)(r-b).$$
The relevant modules   are the corresponding image of
$$  S^{k-1}(S^{b-1}A) \ot \La k A\ot (\La b B)^{\ot k-1} \ot B\ot (
\La b C)^{\ot k-1}\ot C
$$ in
$S^{(k-1)b+1}(A\ot B\ot C)$.
\end{proposition}

The proofs are similar to the proof of Strassen's theorem.

\medskip Another variant is obtained by using products
of $T^{\a}_{\a'}$ with different  $\a$'s and permuted $\a'$'s.
 
\end{example}

\section{Nontriviality of the $(r,s)$-coersive equations}\label{nontrivsect}
We study the image of 
$$\Psi^{r,s}: \Lambda^2(S^sA)\ot S^{r-s}A \ot\wedge^{r}B\ot\wedge^{s}B
\ot\wedge^{s}C\ot\wedge^{r}C)\ra S^{r+s}(A\ot B\ot C).
$$
 Recall that
this may be thought of as first embeding
$ \Lambda^2(S^sA)\ot S^{r-s}A \ot\wedge^{r}B\ot\wedge^{s}B
\ot\wedge^{s}C\ot\wedge^{r}C$ in $(A\ot B\ot C)^{\ot r+s}$ 
according the the recipe in \S\ref{rssect} and
then projecting to the symmetric algebra. We   write the inclusion
and projection as follows:
$$\begin{array}{c}
\Lambda^2(S^sA)\ot S^{r-s}A \ot\wedge^{r}B\ot\wedge^{s}B
\ot\wedge^{s}C\ot\wedge^{r}C \\
 \downarrow \\
%\hookrightarrow
S^sA\ot S^{r-s}A\ot S^sA\ot
\wedge^{s}B\ot\wedge^{r-s}B\ot\wedge^{s}B\ot
\wedge^{s}C\ot\wedge^{r-s}C\ot\wedge^{s}C \\
 | |  \\
S^sA\ot \wedge^{s}B\ot\wedge^{s}C\ot
S^{r-s}A\ot\wedge^{r-s}B\ot \wedge^{r-s}C\ot
S^sA\ot\wedge^{s}B\ot\wedge^{s}C \\
 \downarrow \\
%\hookrightarrow
S^{s}(A\ot B\ot C)\ot S^{r-s}(A\ot B\ot C)\ot
S^{s}(A\ot B\ot C) \\
 \downarrow \\
%\twoheadrightarrow
S^{r+s}(A\ot B\ot C).
\end{array}$$
The first two maps are injective, the last one is surjective but
not injective and the problem is to understand whether its kernel
may contain the subspace we are interested in. For this we need to
understand the above maps in detail, which are made of elementary
maps that we   write down explicitly.

First, we have the map
\begin{eqnarray} \nonumber
\wedge^{r}B & \hookrightarrow & \wedge^{s}B\ot\wedge^{r-s}B
\\ \nonumber
f_1\cwedge f_r & \mapsto & \sum_{I=(i_1<\cdots <i_s)}
\varepsilon(I,{\hat I}) f_{i_1}\cwedge f_{i_s}\ot f_{{\hat
{\i}}_1}\cwedge f_{{\hat {\i}}_{r-s}},
\end{eqnarray}
with the following notation: ${\hat I}=({\hat {\i}}_1<\cdots <
{\hat {\i}}_{r-s})$ is the complementary sequence to $I$ in
$(1,\ldots ,r)$, and $\varepsilon(I,{\hat I})$ is the sign of the
permutation  $(1,\ldots ,r)\mapsto (I,{\hat I})$ (a shuffle).

Second, we have the map
\begin{eqnarray} \nonumber
S^sA\ot \wedge^{s}B\ot\wedge^{s}C & \lra & S^{s}(A\ot
B\ot C)
\\ \nonumber
e^s\ot f_1\cwedge f_s\ot g_1\cwedge g_s  & \mapsto &
\sum_{\s\in\cS_s}\varepsilon (\s)(ef_1g_{\s(1)})\cdots
(ef_sg_{\s(s)}).
\end{eqnarray}
In principle, this information is enough to check if a given
irreducible component of 
$$\Lambda^2(S^sA)\ot S^{r-s}A
\ot\wedge^{r}B\ot\wedge^{s}B \ot\wedge^{s}C\ot\wedge^{r}C
$$
is mapped to zero, or to an isomorphic copy inside $S^{r+s}(A\ot
B\ot C)$. And we just need to test this alternative on some
highest weight vector.

\medskip

Recall the decomposition formulas 
  (e.g. \cite{Macdonald}):
\begin{align}
\nonumber \La 2(S^s V)&= \bigoplus_{j:{\rm  odd}} S_{2s-j,j}V\\
\label{LRrules} \La aV\ot \La bV&=\bigoplus_{\substack{u+v=a+b\\ v\leq\tmin(a,b)}}\Lambda_{u,v}V\\
\nonumber S_{a_1,a_2}V\ot S_bV&= \bigoplus_{\substack{\rho +\s\leq b \\  a_2+\s\leq a_1
}  } S_{a_1+\rho,a_2+\s,b-\rho-\s}V
\end{align}

\medskip

Since we don't have a closed form formula for $\Lambda_{a_1+\rho,a_2+\s,b-\rho-\s}(B\ot C)$,
or more precisely the factors in it of the form
$\Lambda_{u,v}B\ot \Lambda_{u',v'}C$
  we cannot give a closed form formula for all the possible
relevant factors appearing in $S^{r+s}(A\ot B\ot C)$. Even if we did have such
a list, for any given module, we would still have to check that the resulting
map was nonzero before concluding it was present.

We focus on cases that are of length three in $A$  because those of length
two are inherited from $\s_r(\pp 1\times \pp{b-1}\times\pp{c-1})$ which
is treated in \cite{LWsec}.

For example, note that for $s$   odd and $r\ge 2s$,
   $S_{r-s,s,s}A\subset S_{s,s}A\ot S^{r-s}A$   with
multiplicity one. Also $\wedge^{r}B\ot\wedge^{s}B$ contains
$\Lambda_{r,s}B$ with multiplicity one. We prove that the module
$S_{r-s,s,s}A\ot \Lambda_{r,s}B\ot \Lambda_{r,s}C$ is not
mapped to zero in $S^{r+s}(A\ot B\ot C)$ in many cases.

\medskip

We write down a highest weight vector. The tensor product
$f_1\cwedge f_r\ot f_1\cwedge f_s$ gives a highest weight vector
for
  $\Lambda_{r,s}B$
inside $\wedge^{r}B\ot\wedge^{s}B$, where the $f_i$ define a
weight basis of $B$ such that the ordering of
the weights corresponds to the ordering
of the indicies. Similarly  $g_1\cwedge g_r\ot g_1\cwedge
g_s$  gives a highest weight vector for $\Lambda_{r,s}C$. To
find  a highest weight vector for $S_{r-s,s,s}A$ inside
$S^{r-s}A\ot S^{s}A\ot S^{s}A$ we use Young symmetrizers \cite{weyl}.
The symmetrizer  $c_{(r-s,s,s)}$ 
 applied to each of the 
$s$ factors of $A\ot A\ot A$
yields the highest weight vector
\begin{eqnarray} \nonumber
\Theta = & \sum_{\s_1,\ldots, \s_s\in\FS_3}\varepsilon(\s_1)\cdots
\varepsilon(\s_s)e_1^{r-2s}e_{\s_1(1)}\cdots e_{\s_s(1)}\ot
e_{\s_1(2)}\cdots e_{\s_s(2)}\ot e_{\s_1(3)}\cdots e_{\s_s(3)}.
\end{eqnarray}
where $e_1,e_2,e_3$ is an ordered weight basis for $A$ and
$\ep (\s)$ denotes the sign of the permutation $\s$.
Considering the contributions of the six different permutations in
$\cS_3$, we get
\begin{eqnarray} \nonumber
  \Theta
  =\sum_{\a_1+\cdots+\a_6=s}(-1)^{\a_2+\a_4+\a_6}
\binom{s}{\a} & \\ \nonumber
 \hspace*{2cm}
e_1^{r-2s+\a_1+\a_2}e_2^{\a_3+\a_4}e_3^{\a_5+\a_6} &\ot
e_1^{\a_4+\a_5}e_2^{\a_1+\a_6}e_3^{\a_2+\a_3} \ot
e_1^{\a_3+\a_6}e_2^{\a_2+\a_5}e_3^{\a_1+\a_4}.
\end{eqnarray}
where $\binom s\a=\binom s{\a_1}\cdots \binom s{\a_6}$. Now we
take the tensor product of our three highest weight vectors and  
examine  the tensor $\Theta'$ that we get inside $S^{r+s}(A\ot
B\ot C)$.  To show this tensor is nonzero,
we check that the coefficient of
$$(e_1f_1g_1)\cdots
(e_1f_{r-s}g_{r-s})(e_2f_1g_r)\cdots
(e_2f_{s}g_{r+1-s})(e_3f_rg_1)\cdots (e_3f_{r+1-s}g_{s})$$ is
nonzero. The contributions to this monomial in $\Theta '$ is the
sum of the contributions from terms of the form
$$e_1^{r-2s+\a_1+\a_2}e_2^{\a_3+\a_4}e_3^{\a_5+\a_6}f_{{\Hat I}}
g_{{\Hat J}} \ot
e_1^{\a_4+\a_5}e_2^{\a_1+\a_6}e_3^{\a_2+\a_3}f_{1\cdots s}g_J \ot
e_1^{\a_3+\a_6}e_2^{\a_2+\a_5}e_3^{\a_1+\a_4}f_Ig_{1\cdots s},$$
with some coefficient. The first (resp. second, third)    of the
three terms in this product will contribute
to  $\Theta '$ by a
product of terms of the form $(e_if_jg_k)$, where for each given
$i$, the index $k$ describes a set $A_i$ (resp. $B_i$, $C_i$),
with
\begin{eqnarray}
\nonumber A_1 \cup A_2 \cup A_3 &=& {\Hat J} \\
\nonumber B_1 \cup B_2 \cup B_3 &=& J \\
\nonumber C_1 \cup C_2 \cup C_3 &=& \{1,\ldots ,s\}.
\end{eqnarray}
To contribute to our preferred monomial, we also need the conditions
\begin{eqnarray}
\nonumber A_1 \cup B_1 \cup C_1 &=& \{1,\ldots ,r-s\} \\
\nonumber A_2 \cup B_2 \cup C_2 &=& \{1,\ldots ,s\} \\
\nonumber A_3 \cup B_3 \cup C_3 &=& \{r-s+1,\ldots ,r\}.
\end{eqnarray}
Now consider the index $j$ in the different terms $(e_if_jg_k)$.
We need $j=k$ if $k\le r-s$, and $j={\bar k}:=r+1-k$ otherwise.
This leads to one more set of identities,
\begin{eqnarray}
\nonumber A_1 \cup {\bar A}_2 \cup {\bar A}_3 &=& {\Hat I} \\
\nonumber B_1 \cup {\bar B}_2 \cup {\bar B}_3 &=& \{1,\ldots ,s\}, \\
\nonumber C_1 \cup {\bar C}_2 \cup {\bar C}_3 &=& I,
\end{eqnarray}
where ${\bar A}$ denotes the image of $A$ by the map $k\mapsto
{\bar k}$. Note that all these unions are between pairwise disjoint
sets.

The first two relations involving $C_3$   imply that
$C_3=\emptyset$. Since ${\bar B}_2\subset\{1,\ldots ,s\}$, we
deduce that $B_2=\emptyset$, hence $C_1=A_2$. In particular,
$\a_1=\a_4=\a_6=0$. Since also ${\bar B}_1\subset \{1,\ldots
,s\}$, we get that $A_1=A_0\cup \{s+1,\ldots ,r-s\}$ where $A_0$
is the complement to $A_2\cup B_1$ inside $\{1,\ldots ,s\}$.
Comparing $I$ and ${\Hat I}$ we deduce that ${\bar A}_3=B_1$,
hence $A_3={\bar B}_1$ and $B_3={\bar A}_0\cup{\bar A}_2$. In
particular, $I$ and $J$ are determined by $A_0$ and $A_2$. Note
that once we have $I$ and $J$, we can easily compute the signs
$\varepsilon(I,{\Hat I})$ and $\varepsilon(J,{\Hat J})$. The
result is that
$$\varepsilon(I,{\Hat I})\varepsilon(J,{\Hat J})
=(-1)^{(s+\a_2)(r-s+\a_2)}.$$ We deduce that the total
contribution to our monomial is
$$T_{s,r-2s}:=
\sum_{\a_2+\a_3+\a_5=s}(-1)^{\a_2+(s+\a_2)(r-s+\a_2)}\binom{s}{\a}^2
(r-2s+\a_2)!\a_3!\a_5!\a_5!(\a_2+\a_3)!\a_3!(\a_2+\a_5)!.$$
Indeed, for a given $\a$ we have $\binom{s}{\a}$ choices for $A_0,
A_2$, and once these are fixed, the number of permutations sending
the $e_i$'s to the $g_k$ such that $k\in A_i$ is $\# A_1!\# A_2!\#
A_3! =(r-2s+\a_2)!\a_3!\a_5!$, and so on.

So what remains to prove is that $T_{s,r-2s}\neq 0$. Observe that
since $s$ is odd, the product $(s+\a_2)\a_2$ is even. So  
$$ (-1)^{\a_2+(s+\a_2)(r-s+\a_2)}=\bigg\{\begin{array}{ll}
(-1)^{\a_2} & {\rm for}\; r-s \; {\rm even} \\
-1  & {\rm for}\; r-s \; {\rm odd}.
\end{array}
$$
In particular, $T_{s,r-2s}$ is   nonzero for $r$ even. We are not
able to prove all the remaining cases, but we are able to show:

\begin{lemm}
The integer  
$$T_{s,t}=\frac 1{(s!)^2}\sum_{\a+\b+\g=s}(-1)^{\a}\frac{(\a+t)!(\a+\b)!(\a+\g)!}
{\a!\a!}$$ is nonzero for $s,t$ odd and $t\geq s$.
\end{lemm}

\begin{proof} Write $s=2m+1$.
\begin{tiny}
\begin{align*}
(s!)^2T_{s,t}
&= -\sum_{p=0}^m  \sum_{\b=0}^{2m+1-(2p+1)}
\frac{(2p+1+t)!(2p+1+\b)!(2m+1-\b)!}{(2p+1)!(2p+1)!}
 +\sum_{p=0}^m\sum_{\b=0}^{2m+1-2p}
\frac{(2p+t)!(2p+\b)!(2m+1-\b)!}{(2p)!(2p)!}
\\
&= -\sum_{p=0}^m
%\left\{
\{ \sum_{\b=0}^{2m+1-(2p+1)}
\frac{(2p+1+t)(2p+t)!(2p+1+\b)(2p+\b)!(2m+1-\b)!}{(2p+1)(2p)!(2p+1)(2p)!}
\\
& \ \ - \sum_{\b=0}^{2m+1-(2p+1)}
\frac{(2p+t)!(2p+\b)!(2m+1-\b)!}{(2p)!(2p)!}-\frac{(2p+t)!(2m+1)!(2p)!}{(2p)!(2p)!}
%\right\}
\}
\\
&= -\sum_{p=0}^m
%\left\{
\{ \sum_{\b=0}^{2m+1-(2p+1)}
%\left\[
[ \frac{(2p+t)!(2p+\b)!(2m+1-\b)!}{ (2p)! (2p)!}
 (
\frac{(2p+1+t)(2p+1+\b)}{(2p+1)(2p+1)} -1
 )
%\right\]
]  -\frac{(2p+t)!(2m+1)!(2p)!}{(2p)!(2p)!}
%\right\}
\}
\\
&= -\sum_{p=0}^m  \frac{(2p+t)!}{(2p)!}
%\left\{
\{ \sum_{\b=0}^{2m+1-(2p+1)}
%\left\[
[ \frac{ (2p+\b)!(2m+1-\b)!}{ (2p+1)^2 (2p)!}
 (
(2p+1)(t+\b)+t\b
 )
%\right\]
] -  (2m+1)!
%\right\}
\}
\end{align*}
\end{tiny}
In the case $t\geq 2m+1$ this gives the result immediately just by looking
at the $\b=0$ term in the summation and noting it is not  
the only term.  
\end{proof}

We expect $T_{s,t}$ to be always nonzero
when $s,t>1$ and both are odd, but
were unable to prove it.  Note that
it  would be sufficient to prove the case $t=1$ if we could show we always
have the same sign.

\end{document}